\documentclass[12pt,a4paper]{article}
\usepackage{amsmath,amssymb,amsfonts}
\def\Bbb{\mathbb}

\title{\bf Approximation of rational numbers by Dedekind sums}

\author{Kurt Girstmair}
\date{}

\makeatletter
\let\@@maketitle=\maketitle
\def\maketitle{\def\thispagestyle##1{\relax}\@@maketitle}
\makeatother
\newtheorem{theorem}{Theorem}

\def\BE{\begin{equation}}
\def\EE{\end{equation}}
\def\BD{\begin{displaymath}}
\def\ED{\end{displaymath}}
\def\BA{\begin{array}}
\def\EA{\end{array}}
\def\BEA{\begin{eqnarray*}}
\def\EEA{\end{eqnarray*}}
\def\BI{\bibitem}

\def\Z{\Bbb Z}

\def\R{\Bbb R}

\def\phi{\varphi}
\def\EPS{\varepsilon}

\def\CMOD#1#2#3{#1 \equiv #2 \: \mbox{mod}\: #3}

\def\MB{\mbox}
\def\LD{\ldots}

\def\BQ{``}
\def\EQ{'' }

\def\STOP{\hfill$\Box$}

\def\DED{Dedekind }

\begin{document}
\maketitle

\begin{abstract}
\noindent
Given a rational number $x$ and a bound $\EPS$, we exhibit $m,n$ such that $|x-12 s(m,n)|<\EPS$.
Here $s(m,n)$ is the classical \DED sum and the parameters
$m$ and $n$ are completely explicit in terms of $x$ and $\EPS$.
\end{abstract}

\section*{1. Introduction and result}

Let $n$ be a positive integer and $m\in \Z$, $(m,n)=1$. The classical \DED sum $s(m,n)$ is defined by
\BD
   s(m,n)=\sum_{k=1}^n ((k/n))((mk/n))
\ED
where $((\LD))$ is the \BQ sawtooth function\EQ defined by
\BD
  ((t))=\left\{\begin{array}{ll}
                 t-\lfloor t\rfloor-1/2 & \MB{ if } t\in\R\smallsetminus \Z; \\
                 0 & \MB{ if } t\in \Z
               \end{array}\right.
\ED
(see, for instance, \cite[p. 1]{RaGr}). In the present setting it is more
natural to work with
\BD
 S(m,n)=12s(m,n)
\ED instead.
Observe that $S(m+n,n)=S(m,n)$, so one often considers
only arguments $m$ in the range $1\le m\le n$.

It follows from a result in \cite{Hi} that the set
\BD
  \{S(m,n): n>0, 1\le m\le n, (m,n)=1\}
\ED
is dense in the set $\R$ of real numbers. The proof of \cite{Hi} uses continued
fraction expansions of irrational numbers in order to approximate these numbers by \DED sums.
Here we show that each rational number $x$ can be approximated by a \DED sum $S(m,n)$, where
$m$ and $n$ are completely explicit in terms of $x$ and the bound $\EPS$ for the approximation.
Our result reads as follows.

\begin{theorem} 
\label{t1}

Let $\EPS>0$ be given. Let $j,k$ be integers, $0< j\le k$, $(j,k)=1$. Let $m$ be an integer, $m\ge 2/(k\EPS)+1$,
such that
\BE
\label{0}
 \CMOD{mj}1k.
\EE
Let $l$ be a positive integer and put $t=2m+ln-j(m^2+1)$. Further, put $n=k(m^2+1)$.
Then
\BD
\label{0.1}
  S(mt+1,nt)=l-3-j/k+E,
\ED
with $0< E<\EPS$.
\end{theorem} 

Obviously, $l-3-j/k$ runs through all rational numbers $\ge -3$ if $j,k,l$ take the values
permitted by the theorem. Since $S(-mt-1,nt)=-S(mt+1,nt)$ (see \cite[p.26]{RaGr}),
we obtain approximations of all rational numbers. We consider the following

\medskip
\noindent
{\em Example.} Suppose we want to approximate $7/11$ with  $\EPS=1/100$. So $j=4, k=11$, $l=4$.
We have to choose  $m\ge 2/(11\EPS)+1\approx 19.18$ and
(\ref{0}) requires $\CMOD m{3\,}{11}$. Hence we put $m=25$. Now $n=6886$ and $t= 2\cdot 25 +4\cdot 6886-4\cdot 626=25090$.
Further, $mt+1= 627251$ and $nt= 172769740$. We obtain $S(mt+1,nt)=0.6436247\LD=7/11+E$, with $E\approx 73/10000$.

\section*{Proof of Theorem \ref{t1}}

Our basic tool for the proof of Theorem \ref{t1} is

\begin{theorem} 
\label{t2}
Let $n$ be a positive integer and $m\in \Z$, $(m,n)=1$.
Let $m^*\in\Z$ be an inverse of $m$ mod $n$, i.e., $\CMOD{mm^*}1n$, and suppose $m>m^*$ {\rm (}in particular, $m^*$
may be negative{\rm )}. Then we have, with
$t= m-m^*$,
\BE
\label{1}
  S(mt+1,nt)=-3+ \frac 2{nt}+\frac tn.
\EE
\end{theorem} 

\noindent
{\em Proof.}
Observe that
$(mt+1,nt)=1$, since $(mt+1,t)=1$ and $mt+1\equiv m(m-m^*)+1\equiv m^2$ mod $n$.
Let $d$ be a positive integer and $c\in\Z$, $(c,d)=1$. Suppose that $q=md-nc$ is positive.
The three-term relation of Rademacher and Dieter
connects the \DED sums $S(m,n)$ and $S(c,d)$ in the following way:
\BE
\label{2}
 S(m,n)=S(c,d)+S(r,q)+\frac{n^2+d^2+q^2}{ndq}-3
\EE
(see, for instance, \cite[Lemma 1]{Gi3}). Here $r$ is defined as follows:
Let $j, k$ be integers such that
\BE
 \label{3}
  -cj+dk=1.
\EE
Then
\BE
\label{4}
 r=-nk+mj.
\EE
We put $c=m^*$ and $d=n$. By assumption, $t=m-m^*> 0$, so
$q=mn-nm^*=nt>0$. By our choice of $c$ and $d$, we may assume that  $j=-m$ satisfies (\ref{3}) with an appropriate
integer $k$, i.e.
\BD
\label{5}
 m^*m+nk=1.
\ED
If we use this identity to insert for $nk$ in (\ref{4}), we see that
\BD
 r= m^*m-1-m^2=-mt-1.
\ED
It is well-known that $S(m,n)=S(m^*,n)$ (see \cite[p.26]{RaGr}). Therefore, (\ref{2}) reads
\BD
  0=S(-mt-1,nt)+\frac{2n^2+n^2t^2}{n^3t}-3.
\ED
On observing $S(mt+1,nt)=-S(-mt-1,nt)$, we have the desired identity.
\STOP

\medskip
\noindent
{\em Proof of Theorem \ref{t1}.} In the setting of Theorem \ref{t1}, the number $-m+jn/k$ is an integer,
since $n=k(m^2+1)$. Further, it is an inverse of $m$ mod $n$. Indeed, $m(-m+jn/k)=-m^2+njm/k$ and $-m^2=1-n/k$,
so
\BE
\label{6}
  m(-m+jn/k)=1+n(mj-1)/k.
\EE
By (\ref{0}), $mj-1$ is divisible by $k$, so (\ref{6}) says
\BD
  m(-m+jn/k)\equiv 1 \MB{ mod } n.
\ED
In view of Theorem \ref{t2}, we may choose
\BD
   m^*=-m+jn/k-ln.
\ED
Then $t=m-m^*=2m+ln-jn/k=2m+ln-j(m^2+1)$ is positive,
since $j(m^2+1)\le k(m^2+1)=n\le ln.$ Therefore, Theorem \ref{t2} gives
\BD
 S(mt+1,nt)=-3+\frac 2{nt}+\frac tn=-3+\frac 2{nt}+\frac{2m}n+l-\frac jk.
\ED
So this \DED sum equals $l-3-j/k$ up to the error term $2m/n+2/(nt)$.
Now $n>km^2$ and $t\ge 2m$, hence
\BE
\label{7}
  0<\frac{2m}n+\frac 2{nt}<\frac 2{km} +\frac 1{km^3}=\frac{2m^2+1}{km^3}.
\EE
Observe
\BD
   \frac{m^3}{2m^2+1}>\frac{m-1}2\ge \frac 1{k\EPS}.
\ED
Therefore, (\ref{7}) shows $2m/n+2/(nt)<\EPS$.
\STOP


\vspace{0.5cm}
\noindent
Kurt Girstmair            \\
Institut f\"ur Mathematik \\
Universit\"at Innsbruck   \\
Technikerstr. 13/7        \\
A-6020 Innsbruck, Austria \\
Kurt.Girstmair@uibk.ac.at

\end{document}